\documentclass[11pt]{article}

\usepackage{amsmath}
\usepackage{amsthm}
\usepackage{amssymb}
\usepackage{latexsym}
\usepackage{pstricks}
\usepackage{graphicx, pst-plot, pst-node, pst-text, pst-tree}
\usepackage{titlefoot}
\usepackage{enumerate}
\usepackage{gastex}
\usepackage{titlesec}
\usepackage[small,it]{caption}
\usepackage{color}
\definecolor{lightgray}{rgb}{0.9, 0.9, 0.9}
\definecolor{darkgray}{rgb}{0.7, 0.7, 0.7}
\definecolor{darkblue}{rgb}{0, 0, .4}

\usepackage[bookmarks,breaklinks]{hyperref}
\hypersetup{
	colorlinks=true,
	linkcolor=darkblue,
	anchorcolor=darkblue,
	citecolor=darkblue,
	pagecolor=darkblue,
	urlcolor=darkblue,
	pdfpagemode=UseThumbs,
	pdftitle={Simple permutations: decidability and unavoidable substructures},
	pdfsubject={Combinatorics},
	pdfauthor={Robert Brignall, Nik Ruskuc, and Vincent Vatter},
	pdfkeywords={permutation class, restricted permutation, simple permutation}
}

\newtheorem{theorem}{Theorem}[section]
\newtheorem{proposition}[theorem]{Proposition}
\newtheorem{lemma}[theorem]{Lemma}

\newtheorem{conjecture}[theorem]{Conjecture}
\newtheorem{question}[theorem]{Question}
\newenvironment{proofsketch}{\noindent{\it Sketch of proof.}}{\qed\bigskip}

\newcounter{todocounter}

\newcommand{\minisec}[1]{\bigskip\noindent{\bf #1.}}

\makeatletter
\newfont{\footsc}{cmcsc10 at 8truept}
\newfont{\footbf}{cmbx10 at 8truept}
\newfont{\footrm}{cmr10 at 10truept}
\renewenvironment{abstract}%
		{
		  \begin{list}{}%
		     {\setlength{\rightmargin}{1in}%
		      \setlength{\leftmargin}{1in}}%
		   \item[]\ignorespaces\begin{small}}%
		 {\end{small}\unskip\end{list}}

\datefoot{\today}
\amssubj{03B25, 03D05, 05A15, 68F10, 68Q45}

\newpagestyle{main}[\small]{
	\headrule
	\sethead[\usepage][][]
	{\sc Simple Permutations: Decidability and Unavoidable Substructures}{}{\usepage}}

\title{\sc{Simple Permutations: Decidability and Unavoidable Substructures}}
\author{\sc{Robert Brignall, Nik Ru\v{s}kuc, and Vincent Vatter}\thanks{Supported by EPSRC grant GR/S53503/01.}\\
\small School of Mathematics and Statistics\\[-3pt]
\small University of St Andrews\\[-3pt]
\small St Andrews, Fife, Scotland\\[-3pt]
\small \texttt{\{robertb, nik, vince\}@mcs.st-and.ac.uk}\\[-3pt]
\small \texttt{http://turnbull.mcs.st-and.ac.uk/\~{}\{\href{http://turnbull.mcs.st-and.ac.uk/~robertb}{robertb}, \href{http://turnbull.mcs.st-and.ac.uk/~nik}{nik}, \href{http://turnbull.mcs.st-and.ac.uk/~vince}{vince}\}}\\[-10pt]}

\titleformat{\section}
	{\large\sc}
	{\thesection.}{1em}{}	

\date{}

\begin{document}
\maketitle

\pagestyle{main}

\newcommand{\Av}{\operatorname{Av}}
\newcommand{\C}{\mathcal{C}}
\newcommand{\rect}{\operatorname{rect}}

\begin{abstract}
We prove that it is decidable if a finitely based permutation class contains infinitely many simple permutations, and establish an unavoidable substructure result for simple permutations: every sufficiently long simple permutation contains an alternation or oscillation of length $k$.
\end{abstract}

\section{Introduction}\label{sp3-intro}

Simple permutations are the building blocks of permutation classes.  As such, classes with only finitely many simple permutations, e.g., the class of $132$-avoiding permutations, have nice properties.  To name three: these classes have algebraic generating functions (as established by Albert and Atkinson~\cite{albert:simple-permutat:}; see Brignall, Huczynska, and Vatter~\cite{brignall:simple-permutat:} for extensions), are partially well-ordered (see the conclusion), and are finitely based \cite{albert:simple-permutat:}.  It is natural then to ask which finitely based classes contain only finitely many simple permutations.  Our main result establishes that this can be done algorithmically:

\begin{theorem}\label{sp3-main}
It is possible to decide if a permutation class given by a finite basis contains infinitely many simple permutations.
\end{theorem}

\minisec{Permutation classes}
Two sequences $u_1,\dots, u_k$ and $w_1,\dots,w_k$ of distinct real numbers are said to be {\it order isomorphic\/} if they have the same relative comparisons, that is, if $u_i<u_j$ if and only if $w_i<w_j$.
The permutation $\pi$ is said to {\it contain\/} the permutation $\sigma$, written $\sigma\le\pi$, if $\pi$ has a subsequence that is order isomorphic to $\sigma$; otherwise $\pi$ is said to {\it avoid\/} $\sigma$.  For example, $\pi=891367452$ contains $\sigma=51342$, as can be seen by considering the subsequence $91672$ ($=\pi(2),\pi(3),\pi(5),\pi(6),\pi(9)$).  This pattern-containment relation is a partial order on permutations.  We refer to downsets of permutations under this order as {\it permutation classes\/}.  In other words, if $\C$ is a permutation class, $\pi\in\C$, and $\sigma\le\pi$, then $\sigma\in\C$.
We denote by $\C_n$ the set $\C \cap S_n$, i.e.\ the permutations in $\C$ of length
$n$, and we refer to $\sum |\C_n| x^n$ as the {\it generating function for $\C$\/}.
Recall that
an {\it antichain} is a set of pairwise incomparable elements.
For any permutation class $\C$, there is a unique (possibly infinite) antichain $B$ such that $\C=\Av(B)=\{\pi: \beta \not \leq\pi\mbox{ for all } \beta \in B\}$. This antichain $B$, which consists of the minimal permutations not in $\C$, is called the {\it basis} of $\C$.  Permutation classes arise naturally in a variety of settings, ranging from sorting (see, e.g., B\'ona's survey~\cite{bona:a-survey-of-sta:}) to algebraic geometry (see, e.g., Lakshmibai and Sandhya~\cite{lakshmibai:criterion-for-s:}).

It will also be useful to have a pictorial description of order isomorphism.  Two sets $S$ and $T$ of points in the plane are said to be order isomorphic if the axes can be stretched and shrunk in some manner to map one of the sets onto the other, i.e., if there are strictly increasing functions $f,g:\mathbb{R}\rightarrow\mathbb{R}$ such that $\{(f(s_1),g(s_2)) : (s_1,s_2)\in S\}=T$.  (As the inverse of a strictly increasing function is also strictly increasing, this is an equivalence relation).  The {\it plot\/} of the permutation $\pi$ is the point set $\{(i,\pi(i))\}$, and every finite point set in the plane in which no two points share a coordinate (often called a {\it generic\/} or {\it noncorectilinear\/} set) is order isomorphic to the plot of a unique permutation; in practice we simply say that a point set is order isomorphic to a permutation.  This geometric viewpoint makes clear (if they were not already) several symmetries of the pattern-containment order.  The maps $(x,y)\mapsto (-x,y)$ and $(x,y)\mapsto (y,x)$, which when applied to generic point sets correspond to reversing and inverting permutations, generate the dihedral group with eight elements.

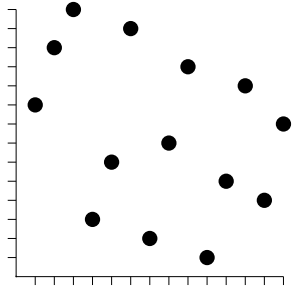
\begin{figure}
\begin{center}
\psset{xunit=0.01in, yunit=0.01in}
\psset{linewidth=0.005in}
\begin{pspicture}(0,0)(140,140)
\psaxes[dy=10,Dy=1,dx=10,Dx=1,tickstyle=bottom,showorigin=false,labels=none](0,0)(140,140)
\pscircle*(10,90){0.04in}
\pscircle*(20,120){0.04in}
\pscircle*(30,140){0.04in}
\pscircle*(40,30){0.04in}
\pscircle*(50,60){0.04in}
\pscircle*(60,130){0.04in}
\pscircle*(70,20){0.04in}
\pscircle*(80,70){0.04in}
\pscircle*(90,110){0.04in}
\pscircle*(100,10){0.04in}
\pscircle*(110,50){0.04in}
\pscircle*(120,100){0.04in}
\pscircle*(130,40){0.04in}
\pscircle*(140,80){0.04in}
\end{pspicture}
\end{center}
\caption{The plot of a simple permutation.}\label{fig-simple-ex}
\end{figure}

\minisec{Simple permutations}
An {\it interval\/} in the permutation $\pi$ is a set of contiguous indices $I=[a,b]$ such that the set of values $\pi(I)=\{\pi(i) : i\in I\}$ also forms an interval of natural numbers.  Every permutation $\pi$ of $[n]=\{1,2,\dots,n\}$ has intervals of size $0$, $1$, and $n$; $\pi$ is said to be {\it simple\/} if it has no other intervals.  Note that simplicity is preserved under the eight symmetries mentioned above.  Figure~\ref{fig-simple-ex} shows the plot of a simple permutation%
\footnote{This assertion suggests a natural question: how quickly can it be decided if a given permutation is simple?  Uno and Yagiura~\cite{uno:fast-algorithms:} present an algorithm that can answer this question in linear time.}.

\begin{figure}
\begin{center}
\psset{xunit=0.02in, yunit=0.02in}
\begin{pspicture}(0,-10)(90,92)
\psframe*[linecolor=lightgray,fillcolor=lightgray,linewidth=0.02in](18,8)(92,82)
\psframe[linecolor=darkgray,linewidth=0.02in](18,8)(92,82)
\psline[linecolor=darkgray,linewidth=0.02in](70,50)(35,50)
\psline[linecolor=darkgray,linewidth=0.02in](60,80)(60,45)
\psline[linecolor=darkgray,linewidth=0.02in](20,70)(65,70)
\psline[linecolor=darkgray,linewidth=0.02in](30,10)(30,75)
\psline[linecolor=darkgray,linewidth=0.02in](90,20)(25,20)
\pscircle*(40,40){0.04in}
\uput[r](40,40){$p_1$}
\pscircle*(50,60){0.04in}
\uput[l](50,60){$p_2$}
\pscircle*(70,50){0.04in}
\uput[r](70,50){$p_3$}
\pscircle*(60,80){0.04in}
\uput[u](60,80){$p_4$}
\pscircle*(20,70){0.04in}
\uput[l](20,70){$p_5$}
\pscircle*(30,10){0.04in}
\uput[d](30,10){$p_6$}
\pscircle*(90,20){0.04in}
\uput[r](90,20){$p_7$}
\pscircle*(80,90){0.04in}
\uput[r](80,90){$x$}
\pscircle*(10,30){0.04in}
\uput[l](10,30){$y$}
\end{pspicture}
\caption{The points $p_1,\dots,p_7$ form a proper pin sequence and the gray box denotes $\rect(p_1,\dots,p_7)$.  The point $x$ satisfies the externality and separation conditions for this pin sequence and thus could be chosen as $p_8$; $y$, however, fails the separation condition.}\label{fig-proper-pin-def}
\end{center}
\end{figure}
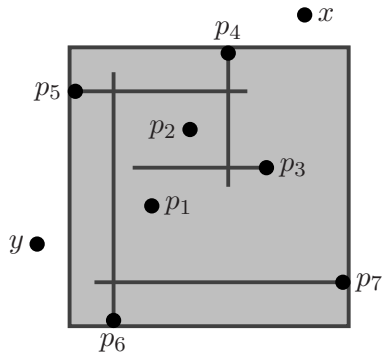

We need several notions from \cite{brignall:simple-permutat:a}.  Given points $p_1,\dots,p_m$ in the plane, we denote by $\rect(p_1,\dots,p_m)$ the smallest axes-parallel rectangle containing them.  A {\it pin\/} for the points $p_1,\dots,p_m$ is any point $p_{m+1}$ not contained in $\rect(p_1,\dots,p_m)$ that lies either horizontally or vertical amongst them.  A {\it proper pin sequence\/} is a sequence of points in the plane satisfying two conditions:\begin{itemize}
\item {\it Externality condition\/}: $p_{i+1}$ must lie outside $\rect(p_1,\dots,p_i)$.
\item {\it Separation condition\/}: $p_{i+1}$ must {\it separate\/} $p_i$ from $\{p_1,\dots,p_{i-1}\}$.  That is, $p_{i+1}$ must lie horizontally or vertically between $\rect(p_1,\dots,p_{i-1})$ and $p_i$.
\end{itemize}
Figure~\ref{fig-proper-pin-def} illustrates these definitions.  The astute reader may note that we have replaced the maximality condition of \cite{brignall:simple-permutat:a} with the externality condition.  This change reflects the differing viewpoints of the papers; while \cite{brignall:simple-permutat:a} was concerned with finding proper pin sequences in permutations, we will be building proper pin sequences from scratch, and in this context the externality and separation conditions together imply the maximality condition.

Proper pin sequences are intimately connected with simple permutations.  In one direction, we have:

\begin{theorem}[Brignall, Huczynska, and Vatter~\cite{brignall:simple-permutat:a}]\label{pins-simple}
If $p_1,\dots,p_m$ is a proper pin sequence then one of the sets of points $\{p_1,\dots,p_m\}$, $\{p_1,\dots,p_m\}\setminus\{p_1\}$, or $\{p_1,\dots,p_m\}\setminus\{p_2\}$ is order isomorphic to a simple permutation.
\end{theorem}

\begin{figure}
\begin{center}
\begin{tabular}{ccccc}
\psset{xunit=0.01in, yunit=0.01in}
\psset{linewidth=0.005in}
\begin{pspicture}(0,0)(120,124)
\psaxes[dy=10,Dy=1,dx=10,Dx=1,tickstyle=bottom,showorigin=false,labels=none](0,0)(120,120)
\pscircle*(10,110){0.04in}
\pscircle*(20,90){0.04in}
\pscircle*(30,70){0.04in}
\pscircle*(40,50){0.04in}
\pscircle*(50,30){0.04in}
\pscircle*(60,10){0.04in}
\pscircle*(70,120){0.04in}
\pscircle*(80,100){0.04in}
\pscircle*(90,80){0.04in}
\pscircle*(100,60){0.04in}
\pscircle*(110,40){0.04in}
\pscircle*(120,20){0.04in}
\end{pspicture}
&\rule{10pt}{0pt}&
\psset{xunit=0.01in, yunit=0.01in}
\psset{linewidth=0.005in}
\begin{pspicture}(0,0)(120,124)
\psaxes[dy=10, Dy=1, dx=10, Dx=1, tickstyle=bottom, showorigin=false,labels=none](0,0)(120,120)
\pscircle*(10,70){0.04in}
\pscircle*(20,50){0.04in}
\pscircle*(30,80){0.04in}
\pscircle*(40,40){0.04in}
\pscircle*(50,90){0.04in}
\pscircle*(60,30){0.04in}
\pscircle*(70,100){0.04in}
\pscircle*(80,20){0.04in}
\pscircle*(90,110){0.04in}
\pscircle*(100,10){0.04in}
\pscircle*(110,120){0.04in}
\pscircle*(120,60){0.04in}
\end{pspicture}
&\rule{10pt}{0pt}&
\psset{xunit=0.01in, yunit=0.01in}
\psset{linewidth=0.005in}
\begin{pspicture}(0,0)(120,124)
\psaxes[dy=10,Dy=1,dx=10,Dx=1,tickstyle=bottom,showorigin=false,labels=none](0,0)(120,120)
\pscircle*(10,20){0.04in}
\pscircle*(20,40){0.04in}
\pscircle*(30,60){0.04in}
\pscircle*(40,80){0.04in}
\pscircle*(50,100){0.04in}
\pscircle*(60,120){0.04in}
\pscircle*(70,90){0.04in}
\pscircle*(80,70){0.04in}
\pscircle*(90,50){0.04in}
\pscircle*(100,30){0.04in}
\pscircle*(110,10){0.04in}
\pscircle*(120,110){0.04in}
\end{pspicture}
\end{tabular}
\end{center}
\caption{From left to right: a parallel alternation, a wedge simple permutation of type $1$, and a wedge simple permutation of type $2$.}\label{fig-par-and-wedge}
\end{figure}
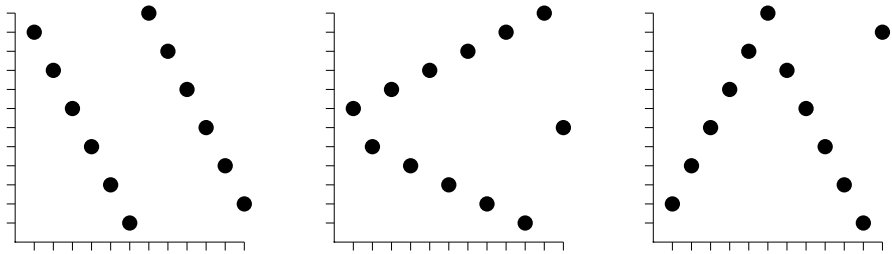

While proper pin sequences are simple or nearly so, there are other types of simple permutations.  We single out three families, plotted in Figure~\ref{fig-par-and-wedge}.  An {\it alternation\/} is a permutation in which every odd entry lies to the left of every even entry, or any symmetry of such a permutation.  A {\it parallel alternation\/} is one in which these two sets of entries form monotone subsequences, either both increasing or both decreasing.  A {\it wedge alternation\/} is one in which the two sets of entries form monotone subsequences pointing in opposite directions.  Whereas every parallel alternation contains a long simple permutation (to form this simple permutation we need, at worst, to remove two points), wedge alternations do not.  However, there are two different ways to add a single point to a wedge alternation to form simple permutations (called {\it simple wedge permutations of types $1$ and $2$\/}).

These families of permutations capture, in a sense formalised below, the diversity of simple permutations.

\begin{theorem}[Brignall, Huczynska, and Vatter~\cite{brignall:simple-permutat:a}]\label{sp2-really-main}
Every sufficiently long simple permutation contains either a proper pin sequence of length at least $k$, a parallel alternation of length at least $k$, or wedge simple permutation of length at least $k$.
\end{theorem}

Theorems~\ref{pins-simple} and \ref{sp2-really-main} show that Theorem~\ref{sp3-main} will follow if we can decide when a class has arbitrarily long parallel alternations, wedge simple permutations, and proper pin sequences.

\section{The Easy Decisions}\label{sp3-easy-decisions}

We begin by describing how to decide if a permutation class given by a finite basis contains arbitrarily long parallel alternations or wedge simple permutations.  Consider first the case of parallel alternations, oriented $\backslash\backslash$, as in Figure~\ref{fig-par-and-wedge}.  These alternations nearly form a chain in the pattern-containment order; precisely, there are two such parallel alternations of each length, and each of these contains a parallel alternation with one fewer point and all shorter parallel alternations of the same orientation.  Thus if the permutation class $\C$ has a basis element contained in any of these parallel alternations, it will contain only finitely many of them.  Conversely, if $\C$ has no such basis element, it will contain all of these alternations.  Therefore we need to characterise the permutations that are contained in any parallel alternation.  To do so, we must first review juxtapositions of classes.

\newcommand{\D}{\mathcal{D}}
\newcommand{\hjuxta}[2]{\left[\begin{array}{cc}#1&#2\end{array}\right]}
\newcommand{\vjuxta}[2]{\left[\begin{array}{c}#1\\#2\end{array}\right]}

Given two permutation classes $\C$ and $\D$, their {\it horizontal juxtaposition\/}, $\left[\begin{array}{cc}\C&\D\end{array}\right]$, consists of all permutations $\pi$ that can be written as a concatenation $\sigma\tau$ where $\sigma$ is order isomorphic to a permutation in $\C$ and $\tau$ is order isomorphic to a permutation in $\D$.

\begin{proposition}[Atkinson~\cite{atkinson:restricted-perm:}]\label{prop-juxta-basis}
Let $\C$ and $\D$ be permutation classes.  The basis elements of the class $\hjuxta{\C}{\D}$ can all be written as concatenations $\rho\sigma\tau$ where either:
\begin{itemize}
\item $\sigma$ is empty, $\rho$ is order isomorphic to a basis element of $\C$, and $\tau$ is order isomorphic to a basis element of $\D$, or
\item $|\sigma|=1$, $\rho \sigma$ is order isomorphic to a basis element of $\C$, and $\sigma\tau$ is order isomorphic to a basis element of $\D$.
\end{itemize}
(In particular, if two classes are finitely based then their juxtaposition is also finitely based.)
\end{proposition}

There is an obvious symmetry to this operation, and the {\it vertical juxtaposition\/} of the classes $\C$ and $\D$ is denoted $\vjuxta{\C}{\D}$.

Proposition~\ref{prop-juxta-basis} is all we need to solve the parallel alternation decision problem.

\begin{proposition}\label{prop-avoid-alt}
The permutation class $\Av(B)$ contains only finitely many parallel alternations if and only if $B$ contains an element of $\Av(123,2413,3412)$ and every symmetry of this class.
\end{proposition}
\begin{proof}
The set of permutations that are contained in at least one (and thus, all but finitely many) parallel alternation(s)  oriented $\backslash\backslash$ is
$$
\hjuxta{\Av(12)}{\Av(12)}=\Av(123,2413,3412),
$$
as desired.
\end{proof}

Like parallel alternations, the wedge simple permutations of a given type and orientation also nearly form a chain in the pattern-containment order, and thus we are able to take much the same approach with them.

\begin{proposition}\label{prop-avoid-wedge1}
The permutation class $\Av(B)$ contains only finitely many wedge simple permutations of type $1$ if and only if $B$ contains an element of 
$$
\Av(1243,1324,1423,1432,2431,3124,4123,4132,4231,4312)
$$
and every symmetry of this class.
\end{proposition}
\begin{proof}
The wedge simple permutations of type $1$ that are oriented $<$, as in Figure~\ref{fig-par-and-wedge}, are contained in
\begin{eqnarray*}
\hjuxta{\vjuxta{\Av(21)}{\Av(12)}}{\{1\}}
&=&
\hjuxta{\Av(132,312)}{\Av(12,21)}
\\
&=&
\Av(1324,1423,1432,2431,3124,4123,4132,4231).
\end{eqnarray*}
It is easy to see that these wedge simple permutations also avoid $1243$ and $4312$, and thus they are contained in the class stated in the proposition, which we call $\D$.

\begin{figure}
\begin{center}
\psset{xunit=0.012in, yunit=0.012in}
\psset{linewidth=0.01in}
\begin{pspicture}(0,-15)(150,150)
% i
\psline[linestyle=dashed](40,0)(40,150)
\psline[linestyle=dashed](40,120)(150,120)
\pscircle*(40,120){0.04in}
% j
\psline[linestyle=dashed](90,0)(90,120)
\psline[linestyle=dashed](40,100)(150,100)
\pscircle*(90,100){0.04in}
% n:
\psline[linestyle=dashed](150,0)(150,150)
\psline[linestyle=dashed](0,70)(150,70)
\pscircle*(150,70){0.04in}
% region labels:
\rput[c](20,35){1}
\rput[c](65,35){2}
\rput[c](120,35){3}
\rput[c](20,110){4}
\rput[c](65,85){5}
\rput[c](65,110){6}
\rput[c](120,85){7}
\rput[c](120,110){8}
\rput[c](95,135){9}
% axes:
\psaxes[dy=1000, Dy=1, dx=1000, Dx=1, tickstyle=bottom, showorigin=false, labels=none](0,0)(150,150)
\psline(40,-3)(40,0)
\psline(90,-3)(90,0)
\psline(150,-3)(150,0)
\rput[c](40,-9){$i$}
\rput[c](90,-9){$j$}
\rput[c](150,-9){$n$}
\end{pspicture}
\end{center}
\caption{The situation in the proof of Proposition~\ref{prop-avoid-wedge1}.}
\label{fig-wedge-simples-basis}
\end{figure}
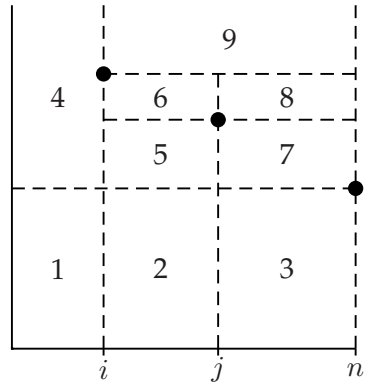

Now take a permutation $\pi\in\D$ of length $n$.  We would like to show that $\pi$ is contained in a wedge simple permutation.  If $\pi\in\vjuxta{\Av(21)}{\Av(12)}$ then $\pi$ is clearly contained in a wedge simple permutation, so suppose this is not the case.  Thus $\pi(1)\cdots\pi(n-1)$ is order isomorphic to a permutation in $\vjuxta{\Av(21)}{\Av(12)}$, and it suffices to show that:
\begin{itemize}
\item the entries of $\pi$ above and to the left of $\pi(n)$ are increasing, and
\item the entries of $\pi$ below and to the left of $\pi(n)$ are decreasing.
\end{itemize}
We prove the first of these items; the second then follows by symmetry because it can be observed from its basis that $\D$ is invariant under complementation, i.e., if the length $n$ permutation $\pi$ lies in $\D$ then so does the permutation $\pi^c$ defined by $\pi^c(i)=n+1-\pi(i)$.  Suppose to the contrary that there is a descent above and to the left of $\pi(n)$.  Thus there are indices $i<j<n$ such that $\pi(i)>\pi(j)>\pi(n)$.  Choose these two indices to be lexicographically minimal with this property.  There must be other entries of $\pi$ as otherwise $\pi$ is simply $321$, which lies in the juxtaposition we have assumed $\pi$ does not lie in.  We now divide $\pi$ into 9 regions as shown in Figure~\ref{fig-wedge-simples-basis}.  About these regions we can state:
\begin{enumerate}
\item[1.] empty because $\pi$ avoids $1432$,
\item[2.] empty because $\pi$ avoids $4132$,
\item[3.] empty because $\pi$ avoids $4312$,
\item[4.] empty because any entry here would necessarily lie below $\pi(j)$ by the minimality of $i$, and would thus create a $2431$ pattern,
\item[5.] empty because $\pi$ avoids $4231$,
\item[6.] empty by the minimality of $j$,
\item[7.] decreasing because $\pi$ avoids $4231$,
\item[8.] empty because $\pi$ avoids $4231$,
\item[9.] increasing because $\pi$ avoids $2431$.
\end{enumerate}
However, this establishes that $\pi$ lies in $\vjuxta{\Av(21)}{\Av(12)}$, a contradiction that completes the proof.
\end{proof}

\begin{proposition}\label{prop-avoid-wedge2}
The permutation class $\Av(B)$ contains only finitely many wedge simple permutations of type $2$ if and only if $B$ contains an element of 
$$
\Av(2134,2143,3124,3142,3241,3412,4123,4132,4231,4312)
$$
and every symmetry of this class.
\end{proposition}
\begin{proofsketch}
The wedge simple permutations of type $2$ that are oriented $\Lambda$, as in Figure~\ref{fig-par-and-wedge}, are contained in
\begin{eqnarray*}
\left[\begin{array}{ccc}
\Av(21)&\Av(12)&\{1\}
\end{array}\right]
&=&
\hjuxta{\Av(213,312)}{\Av(12,21)}\\
&=&
\Av(2134,2143,3124,3142,3241,4123,4132,4231).
\end{eqnarray*}
These wedge simple permutations also avoid $3412$ and $4312$, and an analysis similar to that in the proof of Proposition~\ref{prop-avoid-wedge1} establishes the result.
\end{proofsketch}

\section{Pin Words}\label{sp3-pins}

This leaves only proper pin sequences.  Proper pin sequences, as well as subsets of proper pin sequences, can be described naturally by words over the eight-letter alphabet consisting of the numerals $\{1,2,3,4\}$ and directions $\{L,R,U,D\}$ (standing for left, right, up, and down).  In this section we study these words, laying the groundwork for the proof of our main result in the fourth.

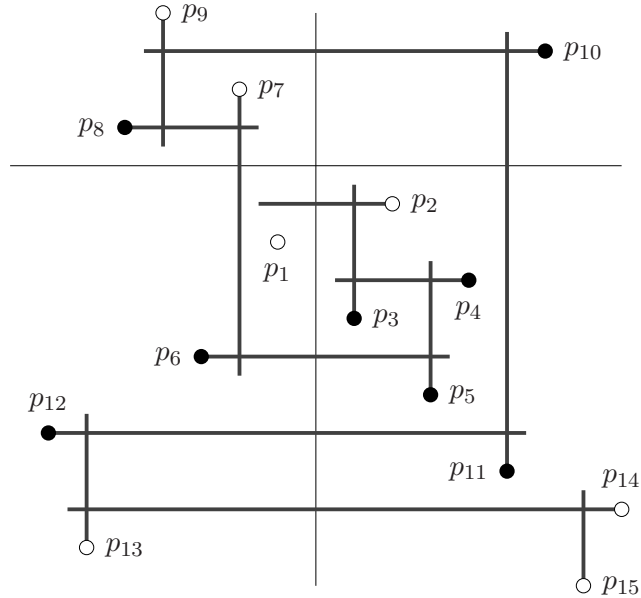
\begin{figure}
\begin{center}
\psset{xunit=0.02in, yunit=0.02in}
\psset{linewidth=0.005in}
\begin{pspicture}(10,7.5)(150,162)
\psaxes[dy=1000,dx=1000,showorigin=false,labels=none](80,120)(160,160)(0,10)
\psline[linecolor=darkgray,linewidth=0.02in](65,110)(100,110)
\psline[linecolor=darkgray,linewidth=0.02in](90,115)(90,80)
\psline[linecolor=darkgray,linewidth=0.02in](85,90)(120,90)
\psline[linecolor=darkgray,linewidth=0.02in](110,95)(110,60)
\psline[linecolor=darkgray,linewidth=0.02in](50,70)(115,70)
\psline[linecolor=darkgray,linewidth=0.02in](60,140)(60,65)
\psline[linecolor=darkgray,linewidth=0.02in](30,130)(65,130)
\psline[linecolor=darkgray,linewidth=0.02in](40,160)(40,125)
\psline[linecolor=darkgray,linewidth=0.02in](140,150)(35,150)
\psline[linecolor=darkgray,linewidth=0.02in](130,40)(130,155)
\psline[linecolor=darkgray,linewidth=0.02in](10,50)(135,50)
\psline[linecolor=darkgray,linewidth=0.02in](20,20)(20,55)
\psline[linecolor=darkgray,linewidth=0.02in](160,30)(15,30)
\psline[linecolor=darkgray,linewidth=0.02in](150,10)(150,35)
\pscircle[fillstyle=solid,fillcolor=white](70,100){0.04in}
\rput[c](70,92){$p_1$}
\pscircle[fillstyle=solid,fillcolor=white](100,110){0.04in}
\rput[l](105,110){$p_2$}
\pscircle*(90,80){0.04in}
\rput[l](95,80){$p_3$}
\pscircle*(120,90){0.04in}
\rput[c](120,82){$p_4$}
\pscircle*(110,60){0.04in}
\rput[l](115,60){$p_5$}
\pscircle*(50,70){0.04in}
\rput[r](45,70){$p_6$}
\pscircle[fillstyle=solid,fillcolor=white](60,140){0.04in}
\rput[l](65,140){$p_7$}
\pscircle*(30,130){0.04in}
\rput[r](25,130){$p_8$}
\pscircle[fillstyle=solid,fillcolor=white](40,160){0.04in}
\rput[l](45,160){$p_9$}
\pscircle*(140,150){0.04in}
\rput[l](145,150){$p_{10}$}
\pscircle*(130,40){0.04in}
\rput[r](125,40){$p_{11}$}
\pscircle*(10,50){0.04in}
\rput[c](10,58){$p_{12}$}
\pscircle[fillstyle=solid,fillcolor=white](20,20){0.04in}
\rput[l](25,20){$p_{13}$}
\pscircle[fillstyle=solid,fillcolor=white](160,30){0.04in}
\rput[c](160,38){$p_{14}$}
\pscircle[fillstyle=solid,fillcolor=white](150,10){0.04in}
\rput[l](155,10){$p_{15}$}
\end{pspicture}
\end{center}
\caption{The proper pin sequence $p_1,\dots,p_{15}$ shown corresponds to the strict pin word $w=3RDRDLULURDLDRD$.  The filled points correspond to the pin word $u=4RDL21DL$, the permutation corresponding to this word, i.e., the permutation order isomorphic to the filled points, is $27453618$.}\label{fig-point-seq}
\end{figure}

The word $w=w_1\cdots w_m\in\{1,2,3,4,L,R,U,D\}^*$ is a {\it pin word\/} if it satisfies:
\begin{enumerate}[(W1)]
\item $w$ begins with a numeral,
\item if $w_{i-1}\in\{L,R\}$ then $w_i\in\{1,2,3,4,U,D\}$, and
\item if $w_{i-1}\in\{U,D\}$ then $w_i\in\{1,2,3,4,L,R\}$.
\end{enumerate}

Pin words with precisely one numeral, which we term {\it strict pin words\/}, correspond to proper pin sequences and it is this correspondence we describe first.  Let $w=w_1\cdots w_m$ denote a strict pin word and begin by placing a point $p_1$ in quadrant $w_1$.  Next take $p_2$ to be a pin in the direction $w_2$ that separates $p_1$ from the origin, denoted $0$.  Continue in this manner, taking $p_{i+1}$ to be a pin in the direction $w_{i+1}$ that satisfies the externality condition and separates $p_i$ from $0,p_1,\dots,p_{i-1}$.  Upon completion, $0,p_1,\dots,p_m$ is a proper pin sequence, and more importantly, $p_1,\dots,p_m$ do as well; it is the latter pin sequence that we say {\it corresponds\/} to $w$.  We say that the {\it permutation corresponding to $w$\/} is the permutation that is order isomorphic to the set of points $p_1,\dots,p_m$.  Conversely, we have the following result.

\begin{lemma}\label{lemma-pin-seq-to-word}
Every proper pin sequence corresponds to a strict pin word.
\end{lemma}
\begin{proof}
Let $p_1,\dots,p_m$ be a proper pin sequence in the plane.  It suffices to place a point $p_0$ (corresponding to the origin) so that $p_0,p_1,\dots,p_m$ form a proper pin sequence.  By symmetry, let us assume that $p_1$ lies below and to the right of $p_2$ and that $p_3$ is a left or right pin.  Hence $p_3$ lies vertically between $p_1$ and $p_2$, and by the separation condition, $p_3$ is the only such pin.  We place $p_0$ vertically between $p_1$ and $p_3$ and minimally to the left of $p_2$, i.e., so that no pin lies horizontally between $p_2$ and $p_0$.  Clearly $p_2$ separates $p_1$ from $p_0$ while $p_3$ separates $p_2$ from $\{p_0,p_1\}$.  Moreover, our placement of $p_0$ guarantees that no later pins separate $\{p_0,p_1,p_2\}$, so since $p_{i+1}$ separates $p_i$ from $\{p_1,\dots,p_{i-1}\}$, it will also separate $p_i$ from $\{p_0,p_1,\dots,p_{i-1}\}$.
\end{proof}

It remains to construct the permutations that correspond to nonstrict pin words.  Letting $w=w_1\cdots w_m$ denote such a word, we begin as before.  Upon reaching a later numeral, say $w_i$, we essentially collapse $p_1,\dots,p_{i-1}$ into the origin and begin anew.  More precisely, we place $p_i$ in quadrant $w_i$ so that it does not separate any of $0,p_1,\dots,p_{i-1}$.  If $w_{i+1}$ is a direction, we take $p_{i+1}$ to be a pin in the direction $w_{i+1}$ that satisfies the externality condition and separates $p_i$ from $0,p_1,\dots,p_{i-1}$; if $w_{i+1}$ is a numeral then we again place $p_{i+1}$ in quadrant $w_{i+1}$ so that it does not separate any of the former points.  In this process we build the {\it sequence of points corresponding to $w$\/}: $p_1,\dots,p_m$.  The {\it permutation corresponding to $w$\/} is again the permutation order isomorphic to this set of points.

We now define an order, $\preceq$, on pin words.  Let $u$ and $w$ be two pin words.  We define a {\it strong numeral-led factor\/} to be a sequence of contiguous letters beginning with a numeral and followed by any number of directions (but no numerals) and begin by writing $u$ in terms of its strong numeral-led factors as $u=u^{(1)}\cdots u^{(j)}$.  We then write $u\preceq w$ if $w$ can be chopped into a sequence of factors $w=v^{(1)}w^{(1)}\cdots v^{(j)}w^{(j)}v^{(j+1)}$ such that for all $i\in[j]$:
\begin{itemize}
\item[(O1)] if $w^{(i)}$ begins with a numeral then $w^{(i)}=u^{(i)}$, and
\item[(O2)] if $w^{(i)}$ begins with a direction, then $v^{(i)}$ is nonempty, the first letter of $w^{(i)}$ corresponds (in the manner described above) to a point lying in the quadrant specified by the first letter of $u^{(i)}$, and all other letters (which must be directions) in $u^{(i)}$ and $w^{(i)}$ agree.
\end{itemize}
Returning a final time to Figure~\ref{fig-point-seq}, the division of $u$ into strong numeral-led factors is $(4RDL)(2)(1DL)$, while $w$ can be written as $(3R)(DRDL)(U)(L)(U)(RDL)(DRD)$.  We now match factors.  Since $w_3$ corresponds to $p_3$ which lies in quadrant $4$, $(4RDL)$ can embed as $(DRDL)$; because $p_8$ lies in quadrant $2$, the $(2)$ factor in $u$ can embed as $(L)$; lastly, $p_{10}$ lies in quadrant $1$, so the $(1DL)$ factor in $u$ can embed as $(RDL)$ in $w$.  This verifies that $u\preceq w$.

%We now define an order $\preceq$ on pin words.  Let $u=u_1\cdots u_k$ and $w=w_1\cdots w_m$ be pin words and suppose that $w$ corresponds to the sequence of points $p_0,p_1,\dots,p_m$ in the manner described above.  Then $u\preceq w$ if there are indices $i_1<\cdots<i_k$ such that:
%\begin{itemize}
%\item[(O1)] if $u_j$ is a direction then $w_{i_j}=u_j$ and $i_{j}=i_{j-1}+1$,
%\item[(O2)] if $u_j$ is a numeral then $p_{i_j}$ lies in quadrant $u_j$, and additionally,
%\item[(O3)] if $u_j$ is a numeral and $w_{i_{j}}$ is a letter then $i_j>i_{j-1}+1$.
%\end{itemize}
%For example, $33L1DL\preceq 1LD23L1RURDL$ (the word from Figure~\ref{fig-point-seq}), as witnessed by the indices $3,5,6,10,11,12$.

This order is not merely a translation of the pattern-containment order on permutations (consider the words $11, 13, 1L, 1D, 21, 23, 2R, 2U,\dots$, which are incomparable under $\preceq$ yet correspond to the same permutation), but $\le$ and $\preceq$ are closely related:

%Note that different pin words can correspond to the same permutation; e.g., the pin words $11, 13, 1L, 1D, 21, 23, 2R, 2U,\dots$ all correspond to the permutation $12$.  Lemma~\ref{pin-words-preceq} addresses this issue.  Before that we have an observation and a definition.
%\begin{lemma}\label{pinseq2}
%Suppose that the pin word $w_1\cdots w_m$ corresponds to the sequence of points $p_0,p_1,\dots,p_m$.  Then for all $i\in[m]$, $p_i$ does not separate any two members of $\{p_0,p_1,\dots,p_{i-2}\}$.
%\end{lemma}
%\begin{proof}
%If $w_i$ is a numeral then this follows from the construction.  If $w_i$ is a direction and if $p_i$ were to separate $\rect(p_0,p_1,\dots,p_{i-2})$ into two parts, then $p_{i-1}$ would lie on one side of this divide, violating the separation condition.
%\end{proof}

\begin{lemma}\label{pin-words-preceq}
If the pin word $w$ corresponds to the permutation $\pi$ and $\sigma\le\pi$ then there is a pin word $u$ corresponding to $\sigma$ with $u\preceq w$.  Conversely, if $u\preceq w$ then the permutation corresponding to $u$ is contained in the permutation corresponding to $w$.
%If $u\preceq w$ then the permutation corresponding to $u$ is contained in the permutation corresponding to $w$.  Conversely, if $\sigma\le\pi$ and $\pi$ corresponds to the pin word $w$ then there is a pin word $u$ corresponding to $\sigma$ with $u\preceq w$.
\end{lemma}
\begin{proof}
If $w=w_1\cdots w_m$ corresponds to the sequence of points $p_1,\dots,p_m$ then the sequence $p_1,\dots,p_{\ell-1},p_{\ell+1},\dots,p_m$ corresponds to the pin word $w_1\cdots w_{\ell-1}w_{\ell+1}'w_{\ell+2}\cdots w_m\preceq w$, where $w_{\ell+1}'$ is the numeral corresponding to the quadrant containing $p_{\ell+1}$.  Iterating this observation proves the first half of the lemma.

The other direction follows similarly.  Write $u$ in terms of its strong numeral-led factors as $u=u^{(1)}\cdots u^{(j)}$ and suppose that the expression $w=v^{(1)}w^{(1)}\cdots v^{(j)}w^{(j)}v^{(j+1)}$ satisfies (O1) and (O2).  Now delete every point in the sequence of points corresponding to $w$ that comes from a letter in a $v^{(i)}$ factor.  By conditions (O1) and (O2) and the remarks in the previous paragraph, it follows that the resulting sequence of points corresponds to $u$.  Therefore the permutation corresponding to $u$ is contained in the permutation corresponding to $w$.
%Now we delete every point corresponding to a letter in a $v^{(i)}$ factor of $w$.  If $w^{(i)}$ begins with a numeral, then (O1) insists that $w^{(i)}=u^{(i)}$, while if $w^{(i)}$ begins with a direction, then (O2) ensures that the resulting sequence of points will correspond to a factor equal to $u^{(i)}$.  Therefore by deleting these points we are reduced to a sequence of points that corresponds to $u$, as desired.
\end{proof}

\section{Brief Review of Regular Languages and Automata}\label{sp3-automata}
\renewcommand{\L}{\mathcal{L}}
\newcommand{\transition}[1]{\stackrel{#1}{\longrightarrow}}

The classic results mentioned here are covered more comprehensively in many texts, for example, Hopcroft, Motwani, and Ullman~\cite{hopcroft:introduction-to:a}, so we give only the barest details.

\begin{figure}
\begin{scriptsize}
\begin{center}
\begin{picture}(52,42)(39,-76)
\gasset{loopCW=n,ELside=r}

% Nodes
\node[Nmarks=i,Nw=12.0,Nh=12.0,Nmr=6.0](start2)(65,-40){}
\node[NLangle=0.0,Nw=12.0,Nh=12.0,Nmr=6.0](c2)(45,-70){{\large $V$}}
\node[NLangle=0.0,Nw=12.0,Nh=12.0,Nmr=6.0](c1)(85,-70){{\large $H$}}

%Edges
\drawedge[ELdist=0.5,syo=1.0,eyo=1.0](c1,c2){$L,R$}
\drawedge[ELdist=1,syo=-1.0,eyo=-1.0](c2,c1){$U,D$}
\drawedge[ELdist=-4](start2,c2){\rotatebox{57}{$1,2,3,4$}}
\drawedge[ELside=l,ELdist=-4](start2,c1){\rotatebox{-57}{$1,2,3,4$}}

\end{picture}
\end{center}
\end{scriptsize}
\caption{The automaton that accepts the language of strict pin words ($V$ and $H$ are accept states).}\label{fig-pin-automaton}
\end{figure}
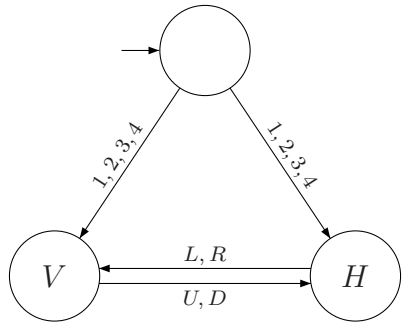

A {\it nondeterministic finite automaton\/} over the alphabet $A$ consists of a set $S$ of {\it states\/}, one of which is designated the {\it initial state\/}, a {\it transition function\/} $\delta$ from $S\times (A\cup\{\varepsilon\})$ into the power set of $S$, and a subset of $S$ designated as {\it accept states\/}.  The {\it transition diagram\/} for this automaton is a directed graph on the vertices $S$, with an arc from $r$ to $s$ labelled by $a$ precisely if $s\in\delta(r,a)$.  The initial state is designated by an inward-pointing arrow.  An automaton {\it accepts\/} the word $w_1\cdots w_m$ if there is a walk from the initial state to an accept state whose arcs are labelled (in order) by $w_1,\dots,w_m$; the set of all such words is the {\it language accepted\/} by the automaton.  For example, Figure~\ref{fig-pin-automaton} shows the transition diagram for the automaton that accepts strict pin words (in this automaton, all states are accept states).

A language that is accepted by a finite automaton is called {\it recognisable\/}.  By Kleene's theorem, the recognisable languages are precisely the {\it regular languages\/}\footnote{The reader unfamiliar with formal languages is welcomed to take this as the definition of regular languages.}, and they have numerous closure properties, of which we use two: the union of two regular languages and the set-theoretic difference of two regular languages are also regular languages.  The other result we need about regular languages is below.

\begin{proposition}\label{regular-finite}
It can be decided whether a regular language given by a finite accepting automaton is infinite.
\end{proposition}
\begin{proofsketch}
A regular language is infinite if and only if one can find a walk in the given accepting automaton that begins at the initial state, contains a directed cycle, and ends at an accept state.
\end{proofsketch}

A {\it finite transducer\/} is a finite automaton that can both read and write.  Transducers also have states, $S$, one of which is designated the initial state and several may be designated accept states.  The transition function for a transducer over the alphabet $A$ is a map from $S\times(A\cup\{\varepsilon\})\times(A\cup\{\varepsilon\})$ into the power set of $A$.  In the transition diagram of a transducer we label arcs by pairs, so the transition $r\transition{a,b}s$ stands for ``read $a$, write $b$''.   Empty inputs and outputs are allowed, both designated by $\varepsilon$, e.g., $r\transition{\varepsilon,b}s$ means ``read nothing, write $b$''.  A word $w\in A^*$ is {\it produced\/} from the word $u\in A^*$ by the transducer $T$ if there is a walk
$$
s_1\transition{u_1,w_1}s_2\transition{u_2,w_2}s_3\cdots\transition{u_m,w_m}s_{m+1}
$$
in the transition diagram of $T$ beginning at the initial state, ending at an accept state, and such that $u=u_1\cdots u_m$ and $w=w_1\cdots w_m$ (note that these $u_i$'s and $w_i$'s are allowed to be $\varepsilon$).  We denote the set of words that the transducer $T$ produces from set of input words $\L$ by $T(\L)$.

\begin{proposition}\label{transducer-regular}
If $\L$ is a regular language and $T$ is a finite transducer then $T(\L)$ is also regular, and a finite accepting automaton for $T(\L)$ can be effectively constructed.
\end{proposition}
\begin{proofsketch}
Let $M$ denote a finite accepting automaton for $\L$.  Suppose that the states of $M$ are $R$ and the states of $T$ are $S$.  The states of an accepting automaton for $T(\L)$ are then $R\times S$, where there is a transition $(r_1,s_1)\transition{b}(r_2,s_2)$ whenever there are transitions $r_1\transition{a}r_2$ and $s_1\transition{a,b}s_2$ in $M$ and $T$, respectively.
\end{proofsketch}

\section{Decidability}\label{sp3-decide}

We are now in a position to prove our main result.  We wish to decide whether the class $\Av(B)$, where $B$ is finite, contains only finitely many simple permutations.  Propositions~\ref{prop-avoid-alt}--\ref{prop-avoid-wedge2} show how to decide if $\Av(B)$ contains arbitrarily long parallel alternations or wedge simple permutations, so by Theorem~\ref{sp2-really-main} it suffices to decide if $\Av(B)$ contains arbitrarily long proper pin sequences.

Consider a permutation $\pi$ that is order isomorphic to a proper pin sequence and thus, by Lemma~\ref{lemma-pin-seq-to-word}, corresponds to at least one strict pin word, say $w$.  If $\pi\not\in\Av(B)$ then $\pi\ge\beta$ for some $\beta\in B$.  By Lemma~\ref{pin-words-preceq}, $\beta$ corresponds to a pin word $u\preceq w$.  Conversely, if $w\succeq u$ for some $u$ corresponding to $\beta\in B$, then Lemma~\ref{pin-words-preceq} shows that $\pi\ge\beta$.  Therefore the set
$$
\{\mbox{strict pin words $w$} : w\succeq u\mbox{ for some $u$ corresponding to a $\beta\in B$}\}
$$
consists of all strict pin words which represent permutations not in $\Av(B)$, so by removing this set from the regular language of all strict pin words we obtain the language of all strict pin words corresponding to permutations in $\Av(B)$.  In the upcoming lemma we prove that for any pin word $u$, the set  $\{\mbox{strict pin words $w$} : w\succeq u\}$ forms a regular language, and thus the language of strict pin words in $\Av(B)$ is regular.  It remains only to check if this language is finite or infinite, which can be determined by Proposition~\ref{regular-finite}.

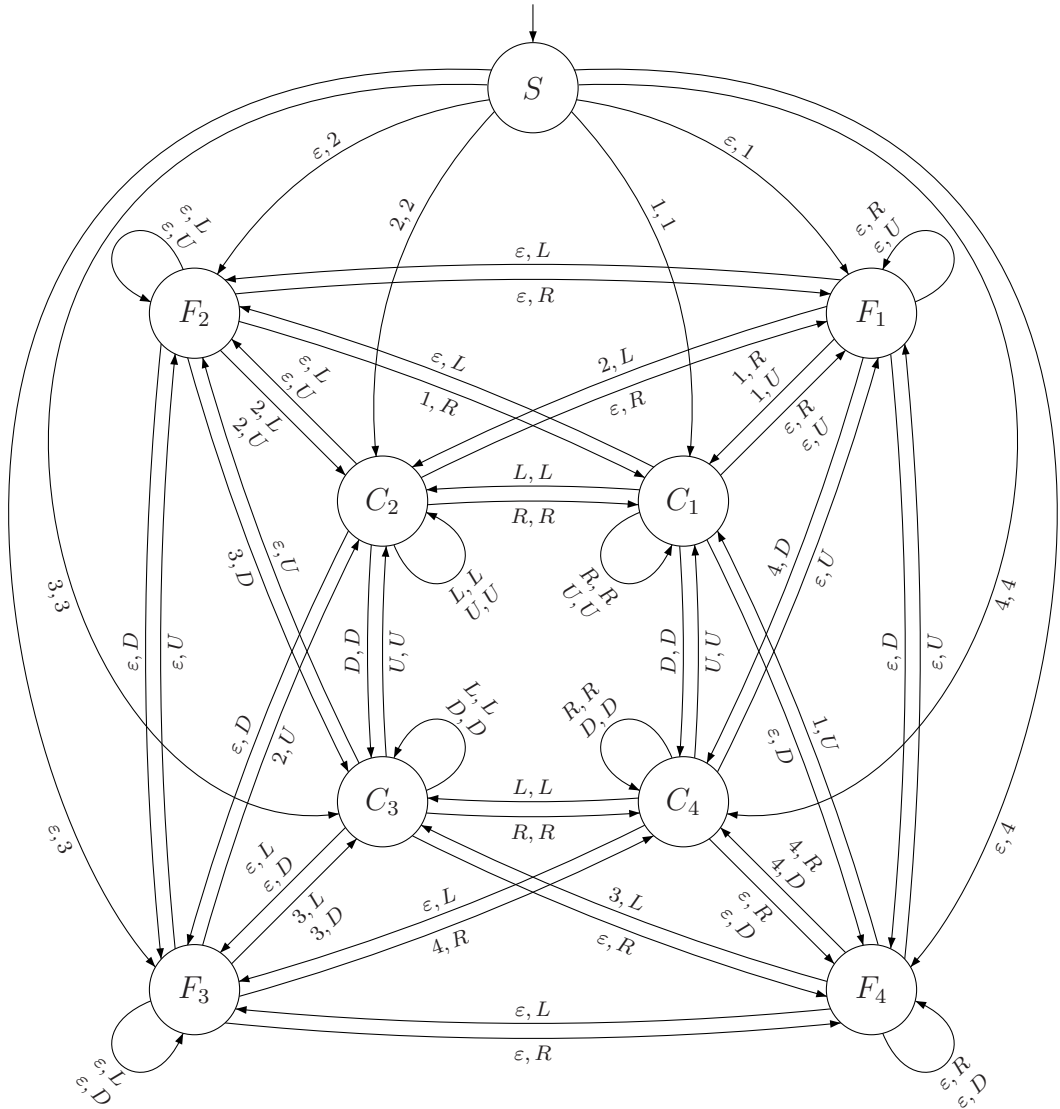
\begin{figure}
\begin{scriptsize}
\begin{center}
\begin{picture}(130,150)(0,-175)\nullfont
\gasset{loopCW=n,ELside=r}
% Nodes
%\node[Nmarks=i,Nw=12.0,Nh=12.0,Nmr=6.0](start)(65,-10){{\large $S_1$}}
\node[Nmarks=i,iangle=90,Nw=12.0,Nh=12.0,Nmr=6.0](start2)(65,-40){{\large $S$}}

\node[NLangle=0.0,Nw=12.0,Nh=12.0,Nmr=6.0](f1)(110,-70){{\large $F_1$}}
\node[NLangle=0.0,Nw=12.0,Nh=12.0,Nmr=6.0](f2)(20,-70){{\large $F_2$}}
\node[NLangle=0.0,Nw=12.0,Nh=12.0,Nmr=6.0](f3)(20,-160){{\large $F_3$}}
\node[NLangle=0.0,Nw=12.0,Nh=12.0,Nmr=6.0](f4)(110,-160){{\large $F_4$}}

\node[NLangle=0.0,Nw=12.0,Nh=12.0,Nmr=6.0](c1)(85,-95){{\large $C_1$}}
\node[NLangle=0.0,Nw=12.0,Nh=12.0,Nmr=6.0](c2)(45,-95){{\large $C_2$}}
\node[NLangle=0.0,Nw=12.0,Nh=12.0,Nmr=6.0](c3)(45,-135){{\large $C_3$}}
\node[NLangle=0.0,Nw=12.0,Nh=12.0,Nmr=6.0](c4)(85,-135){{\large $C_4$}}

% Loops
\drawloop[ELpos=80,ELdist=-4.25,loopangle=45.0](f1){\rotatebox{47}{$\begin{array}{l}\varepsilon,R\\\varepsilon,U\end{array}$}}
\drawloop[ELpos=20,ELdist=-4.25,loopangle=135.0](f2){\rotatebox{-47}{$\begin{array}{l}\varepsilon,L\\\varepsilon,U\end{array}$}}
\drawloop[ELdist=-4.25,loopangle=225.0](f3){\rotatebox{-45}{$\begin{array}{l}\varepsilon,L\\\varepsilon,D\end{array}$}}
\drawloop[ELdist=-4.25,loopangle=-45.0](f4){\rotatebox{45}{$\begin{array}{l}\varepsilon,R\\\varepsilon,D\end{array}$}}

\drawloop[ELdist=-4.25,loopangle=225.0](c1){\rotatebox{-45}{$\begin{array}{l}R,R\\U,U\end{array}$}}
\drawloop[ELdist=-4.25,loopangle=-45.0](c2){\rotatebox{45}{$\begin{array}{l}L,L\\U,U\end{array}$}}
\drawloop[ELdist=-5,loopangle=45.0](c3){\rotatebox{-45}{$\begin{array}{l}L,L\\D,D\end{array}$}}
\drawloop[ELdist=-5,loopangle=135.0](c4){\rotatebox{45}{$\begin{array}{l}R,R\\D,D\end{array}$}}

% Outer edges
\drawedge[ELdist=0.5,syo=4.0,eyo=4.0,curvedepth=-2.5](f1,f2){$\varepsilon,L$}
\drawedge[ELdist=0.5,sxo=-4.0,exo=-4.0,curvedepth=-2.5](f2,f3){\rotatebox{90}{$\varepsilon,D$}}
\drawedge[ELdist=1,syo=-4.0,eyo=-4.0,curvedepth=-2.5](f3,f4){$\varepsilon,R$}
\drawedge[ELdist=1,sxo=4.0,exo=4.0,curvedepth=-2.5](f4,f1){\rotatebox{90}{$\varepsilon,U$}}

\drawedge[ELdist=0.5,sxo=2.0,exo=2.0,curvedepth=2.5](f1,f4){\rotatebox{90}{$\varepsilon,D$}}
\drawedge[ELdist=0.5,syo=-2.0,eyo=-2.0,curvedepth=2.5](f4,f3){$\varepsilon,L$}
\drawedge[ELdist=1,sxo=-2.0,exo=-2.0,curvedepth=2.5](f3,f2){\rotatebox{90}{$\varepsilon,U$}}
\drawedge[ELdist=1,syo=2.0,eyo=2.0,curvedepth=2.5](f2,f1){ $\varepsilon,R$}

% Inner edges
\drawedge[ELdist=0.5,syo=1.0,eyo=1.0,curvedepth=-1](c1,c2){$L,L$}
\drawedge[ELdist=0.5,sxo=-1.0,exo=-1.0,curvedepth=-1](c2,c3){\rotatebox{90}{$D,D$}}
\drawedge[ELdist=1,syo=-1.0,eyo=-1.0,curvedepth=-1](c3,c4){$R,R$}
\drawedge[ELdist=1,sxo=1.0,exo=1.0,curvedepth=-1](c4,c1){\rotatebox{90}{$U,U$}}

\drawedge[ELdist=0.5,sxo=-1.0,exo=-1.0,curvedepth=1](c1,c4){\rotatebox{90}{$D,D$}}
\drawedge[ELdist=0.5,syo=1.0,eyo=1.0,curvedepth=1](c4,c3){$L,L$}
\drawedge[ELdist=1,sxo=1.0,exo=1.0,curvedepth=1](c3,c2){\rotatebox{90}{$U,U$}}
\drawedge[ELdist=1,syo=-1.0,eyo=-1.0,curvedepth=1](c2,c1){$R,R$}

% Edges Outer <-> Inner (matching nodes)
\drawedge[ELdist=-4.25,syo=1.5,eyo=1.5](f1,c1){\rotatebox{45}{$\begin{array}{l}1,R\\1,U\end{array}$}}
\drawedge[ELdist=-4,sxo=-1.5,exo=-1.5](f2,c2){\rotatebox{-45}{$\begin{array}{l}2,L\\2,U\end{array}$}}
\drawedge[ELdist=-4,sxo=1.5,exo=1.5](f3,c3){\rotatebox{45}{$\begin{array}{l}3,L\\3,D\end{array}$}}
\drawedge[ELdist=-4.25,syo=1.5,eyo=1.5](f4,c4){\rotatebox{-45}{$\begin{array}{l}4,R\\4,D\end{array}$}}

\drawedge[ELdist=-4,sxo=1.5,exo=1.5](c1,f1){\rotatebox{45}{$\begin{array}{l}\varepsilon,R\\\varepsilon,U\end{array}$}}
\drawedge[ELdist=-4.25,syo=1.5,eyo=1.5](c2,f2){\rotatebox{-45}{$\begin{array}{l}\varepsilon,L\\\varepsilon,U\end{array}$}}
\drawedge[ELdist=-4.25,syo=1.5,eyo=1.5](c3,f3){\rotatebox{45}{$\begin{array}{l}\varepsilon,L\\\varepsilon,D\end{array}$}}
\drawedge[ELdist=-4,sxo=-1.5,exo=-1.5](c4,f4){\rotatebox{-45}{$\begin{array}{l}\varepsilon,R\\\varepsilon,D\end{array}$}}

% Edges Outer <-> Inner (non-matching nodes)
\drawedge[ELdist=-1.0,syo=2.5,eyo=2.5,curvedepth=-2](f1,c2){\rotatebox{21}{$2,L$}}
\drawedge[ELdist=-0.75,syo=0.5,eyo=0.5,curvedepth=2](c2,f1){\rotatebox{21}{$\varepsilon,R$}}
\drawedge[ELdist=-1.0,sxo=0.5,exo=0.5,curvedepth=2](f1,c4){\rotatebox{69}{$4,D$}}
\drawedge[ELdist=-0.75,sxo=2.5,exo=2.5,curvedepth=-2](c4,f1){\rotatebox{69}{$\varepsilon,U$}}

\drawedge[ELdist=-0.75,syo=0.5,eyo=0.5,curvedepth=2](f2,c1){\rotatebox{-21}{$1,R$}}
\drawedge[ELdist=-1.0,syo=2.5,eyo=2.5,curvedepth=-2](c1,f2){\rotatebox{-21}{$\varepsilon,L$}}
\drawedge[ELdist=-0.75,sxo=-2.5,exo=-2.5,curvedepth=-2](f2,c3){\rotatebox{-69}{$3,D$}}
\drawedge[ELdist=-1.0,sxo=-0.5,exo=-0.5,curvedepth=2](c3,f2){\rotatebox{-69}{$\varepsilon,U$}}

\drawedge[ELdist=-0.75,sxo=-0.5,exo=-0.5,curvedepth=2](f3,c2){\rotatebox{69}{$2,U$}}
\drawedge[ELdist=-1.0,sxo=-2.5,exo=-2.5,curvedepth=-2](c2,f3){\rotatebox{69}{$\varepsilon,D$}}
\drawedge[ELdist=-0.75,syo=-2.5,eyo=-2.5,curvedepth=-2](f3,c4){\rotatebox{21}{$4,R$}}
\drawedge[ELdist=-1.0,syo=-0.5,eyo=-0.5,curvedepth=2](c4,f3){\rotatebox{21}{$\varepsilon,L$}}

\drawedge[ELdist=-1.0,sxo=2.5,exo=2.5,curvedepth=-2](f4,c1){\rotatebox{-69}{$1,U$}}
\drawedge[ELdist=-0.75,sxo=0.5,exo=0.5,curvedepth=2](c1,f4){\rotatebox{-69}{$\varepsilon,D$}}
\drawedge[ELdist=-1.0,syo=-0.5,eyo=-0.5,curvedepth=2](f4,c3){\rotatebox{-21}{$3,L$}}
\drawedge[ELdist=-0.75,syo=-2.5,eyo=-2.5,curvedepth=-2](c3,f4){\rotatebox{-21}{$\varepsilon,R$}}

% Edges from start states
%\drawedge[ELdist=0.5](start,start2){\rotatebox{90}{$1,1$}}

\drawedge[sxo=3,syo=-1,ELpos=35,ELside=l,ELdist=-1.0,curvedepth=7.5](start2,c1){\rotatebox{-62}{$1,1$}}
\drawedge[sxo=-3,syo=-1,ELpos=35,ELdist=-1.0,curvedepth=-7.5](start2,c2){\rotatebox{62}{$2,2$}}

\drawbpedge[ELpos=60,ELside=l,ELdist=-0.5](start2,5,100,c4,-22,48.74){\rotatebox{72}{$4,4$}}
\drawbpedge[ELpos=80,ELside=l,ELdist=-0.75,syo=2,exo=2.0](start2,5,98,f4,38,31.32){\rotatebox{65}{$\varepsilon,4$}}

\drawedge[syo=-1,ELside=l,ELdist=-1.5,curvedepth=7.5](start2,f1){\rotatebox{-34}{$\varepsilon,1$}}
\drawedge[syo=-1,ELdist=-1.5,curvedepth=-7.5](start2,f2){\rotatebox{34}{$\varepsilon,2$}}

\drawbpedge[ELpos=60,ELside=r,ELdist=-0.25](start2,175,100,c3,202,48.74){\rotatebox{-72}{$3,3$}}
\drawbpedge[ELpos=80,ELside=r,ELdist=-0.5,syo=2,exo=-2.0](start2,175,98,f3,142,31.32){\rotatebox{-65}{$\varepsilon,3$}}
\end{picture}
\end{center}
\end{scriptsize}
\caption{The transducer that produces all strict pin words containing the input pin word.}\label{sp3-transducer}
\end{figure}

\begin{lemma}\label{sp3-automaton}
For any pin word $u$, the set $\{\mbox{strict pin words $w$} : w\succeq u\}$ forms a regular language, and a finite accepting automaton for this language can be effectively constructed.
\end{lemma}
\begin{proof}
Let $T$ denote the transducer in Figure~\ref{sp3-transducer}.  We claim that a strict pin word $w$ lies in $T(u)$ if and only if $w\succeq u$.  The lemma then follows by intersecting $T(u)$ with the regular language of all strict pin words.

We begin by noting several prominent features of $T$:
\begin{enumerate}[(T1)]
%\item\label{autobs1} From the initial state $S$, only numerals can be output, so all produced words begin with a numeral.
\item\label{autobs1-5} Every transition writes a symbol.
\item\label{autobs2} Other than the start state $S$, the automaton is divided into two parts, the ``fabrication'' states $F_i$ and the ``copy'' states $C_i$.
\item\label{autobs3} Every transition to a fabrication state has $\varepsilon$ input.
\item\label{autobs4} Every transition from a fabrication state to a copy state reads a numeral and writes a direction, and except for the transitions from $S$, these are the only transitions that read a numeral.
\item\label{autobs5} All transitions between copy states read a direction and write the same direction, these are the only transitions that read a direction, and there is such a transition for every copy state and every direction.
\item\label{autobs6} From every fabrication and copy state, each direction can be output via a transition to a fabrication state with input $\varepsilon$.
\item\label{autobs7} The subscripts of the fabrication and copy states indicate quadrants: if the strict pin word $w_1\cdots w_n$, corresponding to the pin sequence $p_1,\dots,p_n$, has just been written by the transducer and the transducer is currently in state $C_k$ or $F_k$, then $p_n$ lies in quadrant $k$.  Moreover, if the pin word $u_1\cdots u_m$, corresponding to the pin sequence $q_1,\dots,q_m$, has been read and the transducer currently lies in the copy state $C_k$, then $q_m$ lies in quadrant $k$.
\item\label{autobs8} From any state, any copy state can be reached by two transitions, the first being a transition to a fabrication state; for example: $C_2\transition{\varepsilon,D} F_3\transition{4,R} C_4$.
\end{enumerate}

First we prove that $w\succeq u$ for every strict pin word $w$ produced from input $u$ by this transducer.  We prove this by induction on the number of strong numeral-led factors in $u$.  The base case is when $u$ consists of precisely one strong numeral-led factor.  Suppose that the output right before the first letter of $u$ is read is $v^{(1)}$.  There are two cases.  If $v^{(1)}$ is empty, then the transducer is currently in state $S$, and must both read and write the first letter of $u$, moving the transducer into state $C_{u_1}$.  At this point, (T\ref{autobs5}) shows that the transducer could continue to transition between copy states, outputting a word $w=uv^{(2)}\succeq u$.  The only other option available to the transducer (again, by (T\ref{autobs5})) is to transition to a fabrication state, but then (T\ref{autobs4}) shows that the transducer can never again reach a copy state (because $u$ has only one numeral), and thus by (T\ref{autobs3}), it can never finish reading $u$.  In the other case, where $v^{(1)}$ is nonempty, the transducer lies in a fabrication state by (T\ref{autobs4}).  The next transition must then by (T\ref{autobs4}) be into a copy state, and (T\ref{autobs7}) guarantees that the letter written corresponds to a point in quadrant $u_1$.  The same argument as in the previous case shows that the transducer is now confined to copy states until the rest of $u$ has been read, and thus the transducer will output $v^{(1)}w^{(1)}v^{(2)}\succeq u$.

Now suppose that $u$ decomposes into $j\ge 2$ strong numeral-led factors as $u^{(1)}\cdots u^{(j)}$.  By induction, at the point where $u^{(j-1)}$ has just been read, the transducer has output a word $v^{(1)}w^{(1)}\cdots v^{(j-1)}w^{(j-1)}$ and lies in a copy state.  Since the first letter of $u^{(j)}$ is a numeral, the transducer is forced by (T\ref{autobs4}) to transition to a fabrication state, and this transition will write but not read by (T\ref{autobs3}).  The transducer can then transition freely between fabrication states.  Let us suppose that $v^{(1)}w^{(1)}\cdots v^{(j-1)}w^{(j-1)}v^{(j)}$ has been output at the moment just before the transducer begins reading $u^{(j)}$.  As in our second base case above, the transducer must at this point transition to a copy state by (T\ref{autobs4}), which it will do by reading the numeral that begins $u^{(j)}$ and writing a letter that --- by (T\ref{autobs7}) --- corresponds to a point in this quadrant.  The situation is then analogous to the base case, and the transducer will output $v^{(1)}w^{(1)}\cdots v^{(j-1)}w^{(j-1)}v^{(j)}w^{(j)}v^{(j+1)}\succeq u$.

Now we need to verify that the transducer produces every strict pin word $w$ with $w\succeq u$.  Break $u$ into its strong numeral-led factors $u^{(1)}\cdots u^{(j)}$ and suppose that the factorisation $w=v^{(1)}w^{(1)}\cdots v^{(j-1)}w^{(j-1)}v^{(j)}w^{(j)}v^{(j+1)}$ satisfies (O1) and (O2).  If $v^{(1)}$ is nonempty then it can be output immediately by a sequence of transitions to fabrication states by (T\ref{autobs6}); by (O2) and (T\ref{autobs7}), the first letter of $w^{(1)}$ (which must be a direction because $w$ is a strict pin word) can then be output by transitioning to a copy state, from which (T\ref{autobs5}) shows that the rest of $u^{(1)}$ can be read and the rest of $w^{(1)}$ can be written.  If $v^{(1)}$ is empty then $u^{(1)}=w^{(1)}$ by (O1).  The transducer can, by (T\ref{autobs5}), read $u^{(1)}$ and write $w^{(1)}$ by transitioning from $S$ to a copy state and then transitioning between copy states.  Because $w$ is a strict pin word, (O2) shows that $v^{(2)}$ must be nonempty, and (T\ref{autobs6}) shows that $v^{(2)}$ can be output without reading any more letters of $u$.  We then must output $w^{(2)}$ whilst reading $u^{(2)}$.  The only possible obstacle would be reaching the correct copy state, but (T\ref{autobs8}) guarantees that this can be done.  The rest of $u$ can be read, and the rest of $w$ written, in the same fashion.
\end{proof}

The proof of Theorem~\ref{sp3-main} now follows from the discussion at the beginning of the section.

\section{Unavoidable Substructures in Simple Permutations}\label{sp3-unavoidable}

The pin words used to prove our main result allow us to prove an unavoidable substructures result for simple permutations.  This provides an easy-to-check sufficient (but, n.b., not necessary) condition to guarantee that a permutation class contains only finitely many simple permutations.

\begin{figure}
\begin{center}
\psset{xunit=0.01in, yunit=0.01in}
\psset{linewidth=1\psxunit}
\begin{pspicture}(0,0)(150,150)
\psaxes[dy=10, Dy=1, dx=10, Dx=1,tickstyle=bottom,showorigin=false,labels=none]{->}(0,0)(149,149)
\pscircle*(10,40){4.0\psxunit}
\pscircle*(20,10){4.0\psxunit}
\pscircle*(30,60){4.0\psxunit}
\pscircle*(40,30){4.0\psxunit}
\pscircle*(50,80){4.0\psxunit}
\pscircle*(60,50){4.0\psxunit}
\pscircle*(70,100){4.0\psxunit}
\pscircle*(80,70){4.0\psxunit}
\pscircle*(90,120){4.0\psxunit}
\pscircle*(100,90){4.0\psxunit}
\pscircle*(110,140){4.0\psxunit}
\pscircle*(120,110){4.0\psxunit}
\pstextpath[c]{\psline[linecolor=white](125,138)(135,148)}{$\dots$}
\end{pspicture}
\end{center}
\caption{A plot of the increasing oscillating sequence.}
\label{fig-inc-osc}
\end{figure}
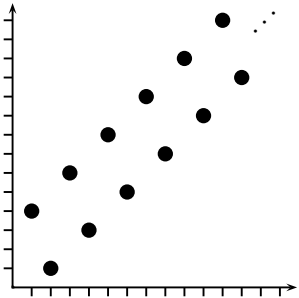

We define the {\it increasing oscillating sequence\/} to be the infinite sequence
$$
4,1,6,3,8,5,\dots,2k+2,2k-1,\dots.
$$
A plot is shown in Figure~\ref{fig-inc-osc}.

We define an {\it increasing oscillation\/} to be any simple permutation that is contained in the increasing oscillating sequence, a {\it decreasing oscillation\/} to be the reverse of an increasing oscillation, and an {\it oscillation\/} to be any permutation that is either an increasing oscillation or a decreasing oscillation.

\begin{theorem}\label{long-pin-seq}
Every sufficiently long simple permutation contains an alternation of length $k$ or an oscillation of length $k$.
\end{theorem}
\begin{proof}
By Theorem~\ref{sp2-really-main}, it suffices to prove that every sufficiently long proper pin sequence contains an alternation or oscillation of length $k$.  Take a proper pin sequence $p_1,\dots,p_m$.  By Lemma~\ref{lemma-pin-seq-to-word}, we may assume that these pins lie in the plane in such a way that $0,p_1,\dots,p_m$ is also a proper pin sequence, where $0$ denote the origin.

We say that this sequence crosses an axis whenever $p_{i+1}$ lies on the other side of the $x$- or $y$-axis from $p_i$, and refer to $\{p_i,p_{i+1}\}$ as a {\it crossing\/}.  First suppose that $p_1,\dots,p_m$ contains at least $2k$ crossings, and so crosses some axis at least $k$ times; suppose that this is the $y$-axis.  Each of these $y$-axis crossings lies either in quadrants $1$ and $2$ or in quadrants $3$ and $4$.  We refer to these as {\it upper crossings\/} and {\it lower crossings\/}, respectively.  By the separation and externality conditions, both pins in an upper crossing lie above all previous crossings, while both pins in a lower crossing lie below all previous crossings.  Thus we can find among the pins of these crossings an alternation of length at least $k$.

Therefore we are done if the pin sequence contains at least $2k$ crossings, so suppose that it does not, and thus that the pin sequence can be divided into at most $2k$ contiguous sets of pins so that each contiguous set lies in the same quadrant.  Each of these contiguous sets is restricted to two types of pin (e.g., a contiguous set in quadrant $3$ can only contain down and left pins) and thus since these two types of pin must alternate, these contiguous sets of pins must be order isomorphic to an oscillation (e.g., a contiguous set in quadrant $3$ must be order isomorphic to an increasing oscillation).  Thus we are also done if one of these contiguous sets has length at least $k$, which it must if the original pin sequence contains at least $m\ge 2k^2$ pins, proving the theorem.
\end{proof}

Thus a class without arbitrarily long alternations or arbitrarily long oscillations necessarily contains only finitely many simple permutations.  First note that these strong conditions are not necessary; for example, the juxtaposition $\hjuxta{\Av(21)}{\Av(12)}$ contains arbitrarily long (wedge) alternations, yet the only simple permutations in this class are $1$, $12$, and $21$.  The work of Albert, Linton, and Ru\v{s}kuc~\cite{insertion} also attests to the strength of these conditions; they prove that classes without long alternations have rational generating functions.

Still, there are benefits to having such a straightforward sufficient condition.  For example, such classes are guaranteed to be partially well-ordered.  As we have already shown how to decide if $\Av(B)$ contains arbitrarily long alternations, to convert Theorem~\ref{long-pin-seq} from a theorem about unavoidable substructures to an easily checked sufficient condition for containing only finitely many simple permutations we need to decide if $\Av(B)$ contains arbitrarily long oscillations.  As with the parallel alternations from Section~\ref{sp3-easy-decisions}, the increasing oscillations nearly form a chain in the pattern-containment order, so we need only compute the class of permutations that are contained in some increasing oscillation, or equivalently, order isomorphic to a subset of the increasing oscillating sequence.  This computation is given without proof in Murphy's thesis~\cite{murphy:restricted-perm:}.  We give a sketch below.

\begin{figure}
\begin{center}
\psset{xunit=0.01in, yunit=0.01in}
\psset{linewidth=1\psxunit}
\begin{pspicture}(0,0)(150,150)
\psaxes[dy=10, Dy=1, dx=10, Dx=1,tickstyle=bottom,showorigin=false,labels=none]{->}(0,0)(149,149)
\pscircle*(10,40){4.0\psxunit}
\pscircle(20,10){4.0\psxunit}
\pscircle(30,60){4.0\psxunit}
\pscircle*(40,30){4.0\psxunit}
\pscircle*(50,80){4.0\psxunit}
\pscircle*(60,50){4.0\psxunit}
\pscircle*(70,100){4.0\psxunit}
\pscircle*(80,70){4.0\psxunit}
\pscircle(90,120){4.0\psxunit}
\pscircle(100,90){4.0\psxunit}
\pscircle(110,140){4.0\psxunit}
\pscircle*(120,110){4.0\psxunit}
\pstextpath[c]{\psline[linecolor=white](125,138)(135,148)}{$\dots$}
\end{pspicture}
\end{center}
\caption{The filled points show the embedding of $2153647$, with rank encoding $1020100$, given by the proof of Proposition~\ref{prop-embed-inc-osc}.}
\label{fig-embed-into-inc-osc}
\end{figure}
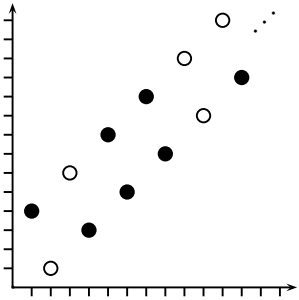

\begin{proposition}\label{prop-embed-inc-osc}
The class of all permutations contained in all but finitely many increasing oscillations is $\Av(321,2341,3412,4123)$.
\end{proposition}
\begin{proofsketch}
It is straightforward to see that every oscillation avoids $321$, $2341$, $3412$, and $4123$, so it suffices to show that every permutation avoiding this quartet is contained in the increasing oscillation sequence.  We use the {\it rank encoding\/}\footnote{We refer the reader to Albert, Atkinson, and Ru\v{s}kuc~\cite{albert:regular-closed-:} for a detailed study of the rank encoding.} for this.  The rank encoding of the permutation $\pi$ of length $n$ is the word $d(\pi)=d_1\cdots d_n$ where
$$
d_i = |\{j : \mbox{$j>i$ and $\pi(j)<\pi(i)$}\}|,
$$
i.e., $d_i$ is the number of points below and to the right of $\pi(i)$.  It is easy to verify that a permutation can be reconstructed from its rank encoding.  Now consider the rank encoding for some $\pi\in\Av(321,2341,3412,4123)$.  It is routine to check that:
\begin{itemize}
\item $d(\pi)\in\{0,1,2\}^*$,
\item $d(\pi)$ does not end in $1$, $2$, or $20$,
\item $d(\pi)$ does not contain $21$, $22$, $111$, $112$, $2011$, or $2012$ factors.
\end{itemize}
We now describe how to embed a permutation with rank encoding satisfying these rules into the increasing oscillating sequence.  Suppose that we have embedded $\pi(1),\dots,\pi(i-1)$.  If $d_i\ge 1$ then we embed $\pi(i)$ as the next even entry in the sequence.  If $d_i=0$ then we embed $\pi(i)$ as the next odd entry if it ends a $20$, $110$, or $2010$ factor, and as the second next odd entry otherwise.  See Figure~\ref{fig-embed-into-inc-osc} for an example.  It remains to show that this is indeed an embedding of $\pi$; to do this it suffices to verify that the number of points of this embedding below and to the right our embedding of $\pi(i)$ is $d_i$.  This follows routinely from the rules above.
\end{proofsketch}

\section{Concluding Remarks}\label{sp3-conclusion}

\minisec{Other contexts}
Analogues of simplicity can be defined for other combinatorial objects, and such analogues have received considerable attention.  For example, let $T$ be a tournament (i.e., directed complete graph) on the vertex set $V(T)$ with (directed) edge set $E(T)$.  For a set $A\subseteq V(T)$ and vertex $v\notin A$, we write $v\rightarrow A$ if $(v,a)\in E(T)$ for all $a\in A$ and similarly $v\leftarrow A$ if $(a,v)\in E(T)$ for all $a\in A$.  An interval in $T$ is a set $A\subseteq V(T)$ such that for all $v\notin A$, either $v\rightarrow A$ or $v\leftarrow A$.  Clearly the empty set, all singletons, and the entire vertex set are all intervals of $T$, and $T$ is said to be simple if it has no others.  Crvenkovi\'c, Dolinka, and Markovi\'c~\cite{crvenkovic:a-survey-of-alg:} survey the algebraic and combinatorial results concerning simple tournaments. 

In the graph case the term ``simple'' is already taken; two correspondent terms are {\it prime\/} and {\it indecomposable\/}.  An {\it interval\/} (also commonly, {\it module\/}) in the graph $G$ is a set $A\subseteq V(G)$ such that every vertex $v\notin A$ is adjacent to every vertex in $A$ or to none.  We refer to Brandst{\"a}dt, Le, and Spinrad's text~\cite{brandstadt:graph-classes:-:} for a survey of simplicity in this context.

Simplicity has also, to some extent, been studied for relational structures in general, for example, by F\"oldes~\cite{foldes:on-intervals-in:} and Schmerl and Trotter~\cite{schmerl:critically-inde:}.

To the best of our knowledge, no analogue of Theorem~\ref{sp3-main} is known for these other contexts.  An approach similar to the one we have taken would require an analogue of Theorem~\ref{sp2-really-main} which, as remarked in \cite{brignall:simple-permutat:a}, remains furtive.

\minisec{Partial well-order}
Recall that a partially ordered set is said to be {\it partially well-ordered\/} ({\it pwo\/}) if it contains neither an infinite strictly decreasing chain nor an infinite antichain.  While permutation classes cannot contain infinite strictly decreasing chains, there are infinite antichains of permutations, see Atkinson, Murphy, and Ru\v{s}kuc~\cite{atkinson:partially-well-:}.  A permutation class with only finitely many simple permutations, on the other hand, is necessarily pwo (Albert and Atkinson~\cite{albert:simple-permutat:} derive this from a result of Higman~\cite{higman:ordering-by-div:}).  Thus Theorem~\ref{sp3-main} bears some resemblance to the pwo decidability question:

\begin{question}
Is it possible to decide if a permutation class given by a finite basis is pwo?
\end{question}

This question is considered in more generality by Cherlin and Latka~\cite{cherlin:minimal-anticha:}. 

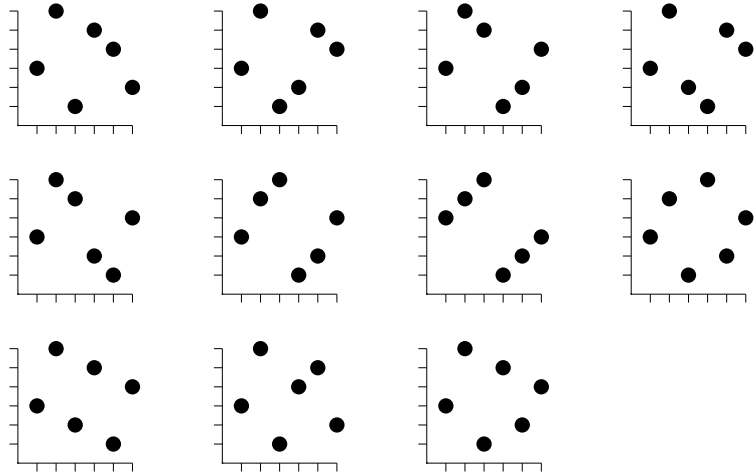
\begin{figure}
\begin{center}
\begin{tabular}{ccccccc}
\psset{xunit=0.01in, yunit=0.01in}
\psset{linewidth=0.005in}
\begin{pspicture}(0,0)(60,64)
\psaxes[dy=10,Dy=1,dx=10,Dx=1,tickstyle=bottom,showorigin=false,labels=none](0,0)(60,60)
\pscircle*(10,30){0.04in}
\pscircle*(20,60){0.04in}
\pscircle*(30,10){0.04in}
\pscircle*(40,50){0.04in}
\pscircle*(50,40){0.04in}
\pscircle*(60,20){0.04in}
\end{pspicture}
&\rule{10pt}{0pt}&
\psset{xunit=0.01in, yunit=0.01in}
\psset{linewidth=0.005in}
\begin{pspicture}(0,0)(60,64)
\psaxes[dy=10,Dy=1,dx=10,Dx=1,tickstyle=bottom,showorigin=false,labels=none](0,0)(60,60)
\pscircle*(10,30){0.04in}
\pscircle*(20,60){0.04in}
\pscircle*(30,10){0.04in}
\pscircle*(40,20){0.04in}
\pscircle*(50,50){0.04in}
\pscircle*(60,40){0.04in}
\end{pspicture}
&\rule{10pt}{0pt}&
\psset{xunit=0.01in, yunit=0.01in}
\psset{linewidth=0.005in}
\begin{pspicture}(0,0)(60,64)
\psaxes[dy=10,Dy=1,dx=10,Dx=1,tickstyle=bottom,showorigin=false,labels=none](0,0)(60,60)
\pscircle*(10,30){0.04in}
\pscircle*(20,60){0.04in}
\pscircle*(30,50){0.04in}
\pscircle*(40,10){0.04in}
\pscircle*(50,20){0.04in}
\pscircle*(60,40){0.04in}
\end{pspicture}
&\rule{10pt}{0pt}&
\psset{xunit=0.01in, yunit=0.01in}
\psset{linewidth=0.005in}
\begin{pspicture}(0,0)(60,64)
\psaxes[dy=10,Dy=1,dx=10,Dx=1,tickstyle=bottom,showorigin=false,labels=none](0,0)(60,60)
\pscircle*(10,30){0.04in}
\pscircle*(20,60){0.04in}
\pscircle*(30,20){0.04in}
\pscircle*(40,10){0.04in}
\pscircle*(50,50){0.04in}
\pscircle*(60,40){0.04in}
\end{pspicture}
\\\\
\psset{xunit=0.01in, yunit=0.01in}
\psset{linewidth=0.005in}
\begin{pspicture}(0,0)(60,64)
\psaxes[dy=10,Dy=1,dx=10,Dx=1,tickstyle=bottom,showorigin=false,labels=none](0,0)(60,60)
\pscircle*(10,30){0.04in}
\pscircle*(20,60){0.04in}
\pscircle*(30,50){0.04in}
\pscircle*(40,20){0.04in}
\pscircle*(50,10){0.04in}
\pscircle*(60,40){0.04in}
\end{pspicture}
&\rule{10pt}{0pt}&
\psset{xunit=0.01in, yunit=0.01in}
\psset{linewidth=0.005in}
\begin{pspicture}(0,0)(60,64)
\psaxes[dy=10,Dy=1,dx=10,Dx=1,tickstyle=bottom,showorigin=false,labels=none](0,0)(60,60)
\pscircle*(10,30){0.04in}
\pscircle*(20,50){0.04in}
\pscircle*(30,60){0.04in}
\pscircle*(40,10){0.04in}
\pscircle*(50,20){0.04in}
\pscircle*(60,40){0.04in}
\end{pspicture}
&\rule{10pt}{0pt}&
\psset{xunit=0.01in, yunit=0.01in}
\psset{linewidth=0.005in}
\begin{pspicture}(0,0)(60,64)
\psaxes[dy=10,Dy=1,dx=10,Dx=1,tickstyle=bottom,showorigin=false,labels=none](0,0)(60,60)
\pscircle*(10,40){0.04in}
\pscircle*(20,50){0.04in}
\pscircle*(30,60){0.04in}
\pscircle*(40,10){0.04in}
\pscircle*(50,20){0.04in}
\pscircle*(60,30){0.04in}
\end{pspicture}
&\rule{10pt}{0pt}&
\psset{xunit=0.01in, yunit=0.01in}
\psset{linewidth=0.005in}
\begin{pspicture}(0,0)(60,64)
\psaxes[dy=10,Dy=1,dx=10,Dx=1,tickstyle=bottom,showorigin=false,labels=none](0,0)(60,60)
\pscircle*(10,30){0.04in}
\pscircle*(20,50){0.04in}
\pscircle*(30,10){0.04in}
\pscircle*(40,60){0.04in}
\pscircle*(50,20){0.04in}
\pscircle*(60,40){0.04in}
\end{pspicture}
\\\\
\psset{xunit=0.01in, yunit=0.01in}
\psset{linewidth=0.005in}
\begin{pspicture}(0,0)(60,64)
\psaxes[dy=10,Dy=1,dx=10,Dx=1,tickstyle=bottom,showorigin=false,labels=none](0,0)(60,60)
\pscircle*(10,30){0.04in}
\pscircle*(20,60){0.04in}
\pscircle*(30,20){0.04in}
\pscircle*(40,50){0.04in}
\pscircle*(50,10){0.04in}
\pscircle*(60,40){0.04in}
\end{pspicture}
&\rule{10pt}{0pt}&
\psset{xunit=0.01in, yunit=0.01in}
\psset{linewidth=0.005in}
\begin{pspicture}(0,0)(60,64)
\psaxes[dy=10,Dy=1,dx=10,Dx=1,tickstyle=bottom,showorigin=false,labels=none](0,0)(60,60)
\pscircle*(10,30){0.04in}
\pscircle*(20,60){0.04in}
\pscircle*(30,10){0.04in}
\pscircle*(40,40){0.04in}
\pscircle*(50,50){0.04in}
\pscircle*(60,20){0.04in}
\end{pspicture}
&\rule{10pt}{0pt}&
\psset{xunit=0.01in, yunit=0.01in}
\psset{linewidth=0.005in}
\begin{pspicture}(0,0)(60,64)
\psaxes[dy=10,Dy=1,dx=10,Dx=1,tickstyle=bottom,showorigin=false,labels=none](0,0)(60,60)
\pscircle*(10,30){0.04in}
\pscircle*(20,60){0.04in}
\pscircle*(30,10){0.04in}
\pscircle*(40,50){0.04in}
\pscircle*(50,20){0.04in}
\pscircle*(60,40){0.04in}
\end{pspicture}
\end{tabular}
\end{center}
\caption{The basis elements of length $6$ for the pin class (up to symmetry).}\label{pin-class-basis-6}
\end{figure}

\minisec{The pin class}
We close with a final, capricious, thought.  The set of permutations that correspond to strict pin words form a permutation class by Lemma~\ref{pin-words-preceq}.  As this class arises from words, it has a distinctly ``regular'' feel, and thus we offer:

\begin{conjecture}\label{conj-pin-class}
The class of permutations corresponding to pin words has a rational generating function.
\end{conjecture}

The enumeration of this class begins $1, 2, 6, 24, 120, 664, 3596, 19004$.  It is not even obvious that this ``pin class'' has a finite basis.  Its shortest basis elements are of length $6$, and there are $56$ of these (see Figure~\ref{pin-class-basis-6}).  The class also has $220$ basis elements of length $7$.

\minisec{Acknowledgements}
We wish to thank Mike Atkinson for fruitful discussions.

\bibliographystyle{acm}
\bibliography{../refs}

\end{document}